\documentclass[12pt,a4paper]{article}
\usepackage{amsmath}
\usepackage{amsfonts}
\usepackage[english, russian]{babel}
\usepackage{amssymb}
\usepackage{amscd}
\usepackage[matrix,arrow,curve]{xy}

\begin{document}
%\begin{flushright}

%\end{flushright}

\begin{center}
\textbf{\Large Cross-frames in Banach spaces} \footnotemark\\[.5cm]
\textbf{S.\,A.~Kreys}\ \ \textrm(Saratov, Russia)\\[.2cm]
\texttt{kreissa@info.sgu.ru}\\[.2cm]
\end{center}
\nopagebreak \footnotetext{
This work was supported by a grant of President of the Russian Federation for young Russian scientists, project MD-300.2011.1, and by the Russian Foundation Basic Research, project 10-01-00097-a.}
\bigskip

%% What about Science School grant SS-4383.2010.1 ? 

\normalsize

This paper generalizes the results obtained in [1] for alternate dual frames in Hilbert spaces on the situation of a Banach space. Additionally, some properties of synthesis operator associated with the alternate dual frame are investigated.

\smallskip

{\bf Definition 1.}
{\hskip 0.2 cm \it A frame for a Hibert space $H$ is a family of vectors
$\left\{\varphi_{n}\right\}_{n=1}^\infty  \subset H\setminus\{0\}$, for which there are positive constants $A$ and $B$, satisfying for all $f \in H$: $$A\left\|f\right\|^2 \leq \sum_{n=1}^\infty
\left|\left(f,\varphi_{n}\right)\right|^2 \leq 
B\left\|f\right\|^2 .$$}

\smallskip

Constants $A$ and $B$ from the above definition are called {\bf\it frame bounds} for a frame, with
$A$ the lower frame bound and $B$ the upper frame bound.
A sequence $\left\{\varphi_{n}\right\}_{n=1}^\infty$ in $H$ with a finite upper frame bound $B$ is called {\it Bessel sequence} in $H$.

\smallskip

{\bf Definition 2.}
{\hskip 0.2 cm \it
A frame $\left\{{y}_{n}\right\}$ in $H$ is called an alternate dual frame for a frame $\left\{{x}_{n}\right\}$,
if for all $x \in H$ it satisfies the following equality:
$$x = \sum_{n=1}^\infty \left(x,y_{n}\right)x_{n}.$$}

\smallskip
In case of a Hilbert space if $\left\{{y}_{n}\right\}$ --- alternate dual frame for $\left\{{x}_{n}\right\}$ then $\left\{{x}_{n}\right\}$ --- alternate dual frame for $\left\{{y}_{n}\right\}$.

\smallskip
The following statements were obtained in [1]:
\smallskip

{\bf Proposition 1.}
{\hskip 0.2 cm \it 
If $\left\{{x}_{n}\right\}$ and $\left\{{y}_{n}\right\}$ are two Bessel sequences in $H$, satisfying $$x = \sum_{n=1}^\infty \left(x,y_{n}\right)x_{n}$$
for any $x \in H$, then $\left\{{x}_{n}\right\}$ and $\left\{{y}_{n}\right\}$ --- alternate dual frames.
}

\smallskip

Given two sequences $\left\{{x}_{n}\right\}$ and $\left\{{y}_{n}\right\}$ in Hilbert space $H$, $\left\{{x}_{n}\right\}$ is called equivalent to $\left\{{y}_{n}\right\}$ if $Tx_n=y_n$ is an invertible linear bounded operator on $H$.

\smallskip

{\bf Proposition 2.}
{\hskip 0.2 cm \it 
Let consider $\left\{{x}_{n}\right\}$ and $\left\{{y}_{n}\right\}$ be the frames in $H$. The family $\left\{{y}_{n}\right\}$ is equivalent to some alternate dual frame for $\left\{{x}_{n}\right\}$ if and only if the operator
$$Sx = \sum_{n=1}^\infty \left(x,x_{n}\right)y_{n}$$ is invertible as an operator from $H$ to $H$. 
}

\smallskip

{\bf Proposition 3.}
{\hskip 0.2 cm \it 
Let $\left\{{x}_{n}\right\}$ be a frame in $H$ and $\left\{{y}_{n}\right\}$ be an alternate dual frame. If $S: H \mapsto H$ --- invertible operator then $\left\{{Sx}_{n}\right\}$ and $\Bigl\{\left(S^{-1}\right)^{\ast}y_{n}\Bigr\}$ are alternate dual frames in $H$.
}

\smallskip

Grochenig [2] first generalized frames onto Banach spaces. There are several different approaches to the definition of a frame in the Banach space.
For example, in [1, Def. 3.3], the frame is defined as a projection of an unconditional basis of the larger Banach space.

\smallskip

Let us consider Banach space $X$ along with its dual space $X^{\ast}$. 
Futhermore, let $X_d$ be a Banach space, consisting of sequences of numbers
$a=\{a_n\}$ and satisfiying the following condition: family of vectors $ \{ \varepsilon_n \} $ forms a basis in $ X_d $, where $ \varepsilon_n = \{ \delta_{nj} \} $ with $\delta_{nj}$ be a Kronecker delta. This basis is often called natural basis. Taking that into account we can consider the dual space $X_d^{\ast}$ again as the space of sequences.\\ One more important requirement is that the dual sequence $ \{ \varepsilon_n^{\ast} \} $ must be a basis in $ X_d^{\ast} $.

\smallskip
The following definition was first introduced in [3] as a notion of $X_d$-frame:
\smallskip

{\bf Definition 3.}
{\hskip 0.2 cm \it 
Consider $\{y_n\} \subset X^{\ast}$. Let us suppose that there exist such positive numbers $A, B$, that the equality 
$$A\left\|x\right\|_{X} \leq \left\|\Bigl\{\left(x,y_{n}\right)\Bigr\}\right\|_{X_{d}} \leq B\left\|x\right\|_{X} \eqno (1)$$
holds for any $x\in X$.
In this case the sequence $\left\{{y}_{n}\right\}$ is called a frame.
}

\vskip 0.4 cm

{\bf Definition 4.}
{\hskip 0.2 cm \it 
Consider $\{x_n\} \subset X$, for any $y\in X^{\ast}$ we have $\left\{\left(x_n, y\right)\right\} \in X_d^{\ast}$, and there exist such positive numbers $\widetilde{A}, \widetilde{B}$, that the equality 
$$\widetilde{A}\left\|y\right\|_{X^{\ast}} \leq \left\|\Bigl\{\left(x_{n},y\right)\Bigr\}\right\|_{X^{\ast}_{d}} \leq \widetilde{B}\left\|y\right\|_{X^{\ast}} \eqno (2)$$
holds for any $y\in X^{\ast}$.
In this case we call the sequence $\left\{{x}_{n}\right\}$ a co-frame.
}

\vskip 0.4 cm

{\bf Definition 5.}
{\hskip 0.2 cm \it 
Suppose we are given $\{x_n\}$ --- co-frame and $\{y_n\}$ --- frame. \\
If for all $x \in X$ and for all $y \in X^{\ast}$ the following formulas recovery
$$x = \sum_{n=1}^\infty \left(x,y_{n}\right)x_{n}, \eqno (3)$$
$$y = \sum_{n=1}^\infty \left(x_{n},y\right)y_{n}, \eqno (4)$$
then we say that a pair of sequences $\left\{x_n\right\}$ and $\left\{y_n\right\}$ form a cross-frame.
}

\vskip 0.4 cm

{\bf Remark.}
{\hskip 0.2 cm  
Cross-frame is essentially an analogue of alternate
dual frames in Hilbert spaces.
Indeed, if we consider (1) and (2) in the case when $X$ --- Hilbert space, then the sequences
$\left\{x_{n}\right\}_{n=1}^\infty$ and $\left\{y_{n}\right\}_{n=1}^\infty$ are frames in $X$ in the sense of Definition 1, and from (3) it follows that they comply with the Definition 2, ie, form a pair of alternate dual frames.
 Therefore, further in the case when sequences $\left\{{x}_{n}\right\}$ and $\left\{{y}_{n}\right\}$ form a cross-frame, sequence $\left\{{x}_{n}\right\}$ will be called alternate dual frame for frame $\left\{{y}_{n}\right\}$, and $\left\{{y}_{n}\right\}$ --- the alternate dual frame for 
 $\left\{{x}_{n}\right\}$. 
}

\smallskip

Let us introduce the operators of synthesis and analysis associated with cross-frame $\Bigl(\left\{{x}_{n}\right\}, \left\{{y}_{n}\right\}\Bigr)$.
Operator $\widetilde{R}: X \mapsto X_d$ defined by the equalities 
$\widetilde{R}x=\Bigl\{\left(x,y_n\right)\Bigr\}$, is called the
analysis of the sequence $\left\{{y}_{n}\right\}$. Note that inequality (1) implies the boundedness of operator $\widetilde{R}$.\\
Let us find the adjoint operator for $\widetilde{R}$. Let $C_{00}$ be the space of finite sequences converging to zero, $b \in C_{00}$:
$$\left(x,\widetilde{R}^{\ast}b\right) = \left(\widetilde{R}x,b\right) = \sum\left(x,y_n\right)b_n = \left(x, \sum{{b_n}{y_n}}\right).$$
Thus, $\widetilde{R}^{\ast}b = \sum{{b_n}{y_n}}$ on the whole $C_{00}$. But since the sequence $\{ {\varepsilon_n}^{\ast} \}$ is a basis $X_d^\ast$ and $C_{00}$ is dense in $X_d^{\ast}$ the operator $\widetilde{R}^{\ast}$ has the following form:  
$$\widetilde{R}^{\ast}b = \widetilde{R}^{\ast}\left({\sum{{b_n}{\varepsilon_n}^{\ast}}}\right)  
= \sum{b_n}\widetilde{R}^{\ast}{\varepsilon_n}^{\ast} = \sum{{b_n}{y_n}}.$$
on the whole space $X_d^\ast$.
Operator $\widetilde{S} = \widetilde{R}^{\ast}$ is called the synthesis of the sequence $\left\{{y}_{n}\right\}$. 
Analysis operator of the sequence $\left\{{x}_{n}\right\}$, acting from $X^{\ast}$ to ${X_d}^{\ast}$, is defined as follows:\\
$Ry=\Bigl\{\left(x_n,y\right)\Bigr\}$. Its boundedness follows from (2).\\ 
Let us define the synthesis operator of the frame $\left\{{x}_{n}\right\}$ as follows:
$$Sa = \sum{{a_n}{x_n}}$$ on all sequences $\{a_n\} \in X_d$, for which the series $\sum{{a_n}{x_n}}$ converges.
Consider the operator $S_0: C_{00} \mapsto X$, where ${S_0}a = Sa$ for $a \in C_{00}$. The chain below shows its boundedness:
$$\left({S_0}a,y\right) = \left(\sum{{a_n}{x_n}},y\right) = \sum{a_n}\left(x_n,y\right) = \left(a,\Bigl\{\left(x_n,y\right)\Bigr\}\right) = \left(a,Ry\right) \eqno (5)$$ with $a \in C_{00}$. 
Due to the fact that the sequence $\{\varepsilon_n\}$
forms a basis in $X_d$ the space $C_{00}$ is dense in $X_d$. This means that there exists unique bounded linear operator from $X_d$ to $X$ which is an extension of $S_0$ on all $X_d$, it is the synthesis operator $S$:
$$\sum{a_n}{x_n} = \sum{a_n}{S_0}{\varepsilon_n} =  \sum{a_n}S{\varepsilon_n} = S{\sum{{a_n}{\varepsilon_n}}} = Sa.$$
And finally from (5) it follows that $S^{\ast} = R$.

\smallskip
We will use these notations for the identity operators in $X$ and $X^{\ast}$ correspondingly:
$I_X$ and $I_{X^{\ast}}$.
\smallskip

{\bf Proposition 4.}
{\hskip 0.2 cm \it 
Formula (4) is not required in the definition of a cross-frame, namely: it can be obtained from (1) --- (3).
}
\\  
{Proof. \hskip 0.2 cm 
The next formula immediately follows from (3):
$$S\widetilde{R} = I_{X}. \eqno (6)$$
Turning to the conjugation in (6), we obtain the equality for adjoint operators:
$${\widetilde{R}}^{\ast}S^{\ast} = I_{X^{\ast}}. \eqno (7)$$
And then due to $${\widetilde{R}}^{\ast}b=\sum_{n=1}^\infty {b_n}{y_n}$$ we obtain the validity of the formula (4) for all $y$ from $X^{\ast}$.
}

\smallskip

We formulate a statement similar to Proposition 1 for cross-frames in Banach spaces.

\smallskip

{\bf Theorem 5.}
{\hskip 0.2 cm \it 
Let $\left\{x_n\right\} \subset X$, $\left\{y_n\right\} \subset X^{\ast}$ and
for all $y$ from $X^{\ast}$, for all $x$ from $X$ we have
$$\Bigl\{\left(x_{n},y\right)\Bigr\} \in X^{\ast}_{d} \eqno (8)$$
and
$$\Bigl\{\left(x,y_{n}\right)\Bigr\} \in X_{d}.$$
 Suppose that we also have the reconstruction formula (3)
$$x = \sum_{n=1}^\infty \left(x,y_n\right)x_n.$$
Then sequences $\left\{x_n\right\}$ and $\left\{y_n\right\}$ form a cross-frame.
}
\\ Proof. 
{\hskip 0.2 cm 
Let $ L: X_d \mapsto X $ be the linear bounded operator acting by the formulas
$L\varepsilon_n = x_n$ with $n \in N$. Since $\left(\varepsilon_n\right)$ is a basis in $X_d$, 
then for any sequence $a \in X_d$ : $a = \sum{a_n}{\varepsilon_n}$, $a = \left(a_n\right)$. 
Lets find the adjoint operator $L^{\ast}$:
$$\left(La, y\right) = \left(L\sum{a_n}{\varepsilon_n}, y\right) = \left(\sum{a_n}{L\varepsilon_n}, y\right)=$$
$$\sum{a_n}\left(x_n, y\right) = \Bigl(\left\{a_n\right\}, \{\left(x_n, y\right)\}\Bigr) = \left(a, L^{\ast}y\right).$$
Thus, 
$$L^{\ast}y = \Bigl\{\left(x_n,y\right)\Bigr\}. \eqno (9)$$
$L^{\ast}$ is acting onto $X_d^{\ast}$ as seen with (8). 
Hence the operator $L^{\ast}$ is also bounded.
	Next, consider the reconstruction formula (3):
$$x = \sum_{n=1}^\infty \left(x,y_{n}\right)x_{n} = \sum_{n=1}^\infty \left(x,y_{n}\right)L\varepsilon_{n} =
L\left(\sum_{n=1}^\infty \left(x,y_{n}\right)\varepsilon_{n}\right) = La, a \in X_d.$$
The last equation means that $L$ - operator on, ie, a surjection.
Then the operator $L^{\ast}$ is injective. Using this fact together with formula (9), we obtain the existence of such constant $A > 0$ that $A\left\|y\right\|_{X^{\ast}} \leq \left\|\Bigl(\left(x_n,y\right)\Bigr)\right\|_{X^{\ast}_{d}}$ for all $y$ from $X^{\ast}$. 
\vskip 0.01 cm
\hskip 0.2 cm 
Let us prove that for a sequence $\left(y_n\right)$ the inequality (1) holds.
By completely repeating the proof of Theorem 4, we can obtain formula (4).
Consider the analysis operator $\widetilde{R}x=\Bigl(\left(x,y_n\right)\Bigr)$. 
Define $B = \left|\left|\widetilde{R}\right|\right| < \infty$, then 
$$\left|\left|\widetilde{R}x\right|\right|_{X_d} = \left|\left|\Bigl(\left(x,y_n\right)_{n=1}^\infty\Bigr)\right|\right|_{X_d} 
\le B\left|\left|x\right|\right|_X.$$
%\vskip 0.01 cm
%\hskip 0.2 cm 
Formula (4) implies that for every $y \in X^{\ast}$ there exists a sequence
$b \in {X_d}^{\ast}$ such that $y = \widetilde{R}^{\ast}b$. 
This means that the operator $\widetilde{R}^{\ast}$ is surjective, hence the operator $\widetilde{R}$ is injective. 
Thus, there exist positive constants $A$ and $B$, that for the sequence $\left(y_n \right)$ we have the bilateral inequality (1).
}

\smallskip

The following theorems generalize Proposition 2 to the situation of Banach spaces.

\smallskip

{\bf Theorem 6.}
{\hskip 0.2 cm \it 
Consider $\left\{x_n\right\} \subset X$ --- co-frame, $\left\{y_n\right\}$ --- frame. 
Then the sequence $\left\{x_n\right\}$ is equivalent to an alternate dual frame for $\left\{y_n\right\}$ if and
only if the operator $$Ux = \sum_{n=1}^\infty \left(x,y_{n}\right)x_{n}$$ is invertible as an operator
from $X$ to $X$.
}
\\ Proof. 
{
\hskip 0.2 cm 
Necessity. Choose an invertible operator $T: X \mapsto X$ such that the sequence
 $\left(Tx_n\right)$ is alternate dual frame for $\left(y_n\right)$. 
 Then for all $x$ from $X$ we will have the representation
$$x = \sum_{n=1}^\infty \left(x,y_{n}\right)Tx_{n}.$$
Hence $U = T^{-1}$ --- invertible operator on $X$.

\vskip 0.01 cm
\hskip 0.2 cm 
Sufficiency. Given invertible operator $U: X \mapsto X$ such that for all $x$ from $X$
$$Ux = \sum_{n=1}^\infty \left(x,y_{n}\right)x_{n},$$
ie
$$x = \sum_{n=1}^\infty \left(x,y_{n}\right)U^{-1}x_{n}.$$
If we put $y = {\left(U^{-1}\right)}^{\ast}y$ in inequality (2), we obtain
$$\widetilde{A}\left\| {\left(U^{-1}\right)}^{\ast}y \right\|_{X^{\ast}} \leq \left\| \Bigl(\left(x_{n}, {\left(U^{-1}\right)}^{\ast}y \right)\Bigr
) \right\|_{X^{\ast}_{d}} \leq \widetilde{B}\left\|  {\left(U^{-1}\right)}^{\ast}y  \right\|_{X^{\ast}}. \eqno (10)$$
Since the operator $U^{-1}$ is invertible and bounded, (10) yields
$$\widetilde{\widetilde{A}}\left\| y \right\|_{X^{\ast}} \leq \widetilde{A}\left\| {\left(U^{-1}\right)}^{\ast}y \right\|_{X^{\ast}} \leq
\left\| \Bigl(\left(x_{n}, {\left(U^{-1}\right)}^{\ast}y \right)\Bigr)\right\|_{X^{\ast}_{d}} \leq \widetilde{B}\left\| {\left(U^{-1}\right)}^{\ast}y \right\|_{X^{\ast}}
\leq \widetilde{\widetilde{B}}\left\| y \right\|_{X^{\ast}}.$$
And with the formula $\left(U^{-1}x_n, y\right) = \left(x_{n}, {\left(U^{-1}\right)}^{\ast}y\right)$ we obtain inequality (2) for the sequence $\left(U^{-1}x_n\right)$.
But then $\left(U^{-1}x_n\right)$ --- an alternate dual frame for the $\left(y_n\right)$.
}

\smallskip

{\bf Theorem 7.}
{\hskip 0.2 cm \it 
Consider $\left\{x_n\right\} \subset X$ --- co-frame, $\left\{y_n\right\}$ --- frame. 
Then the sequence $\left\{y_n\right\}$ is equivalent to an alternate dual frame for $\left\{x_n\right\}$ if and
only if the operator $$Vy = \sum_{n=1}^\infty \left(x_{n},y\right)y_{n}$$ is invertible as an operator from $X^{\ast}$ to $X^{\ast}$.
}

{
\hskip -0.6 cm 
Proof. 
\hskip 0.2 cm 
Necessity. Choose an invertible operator $T: X^{\ast} \mapsto X^{\ast}$ such that the sequence 
 $\left(Ty_n\right)$ is alternate dual frame for $\left(x_n\right)$. Then for all $y$ from $X^{\ast}$ 
 we will have the representation formula (4):
$$y = \sum_{n=1}^\infty \left(x_n,y\right)Ty_{n}.$$
Thus, $V = T^{-1}$ --- an invertible operator on $X^{\ast}$.
\vskip 0.01 cm
\hskip 0.2 cm 
Sufficiency. 
Given invertible operator $V: X^{\ast} \mapsto X^{\ast}$ such that for all $y$ from $X^{\ast}$
$$Vy = \sum_{n=1}^\infty \left(x_n,y\right)y_{n},$$
ie
$$y = \sum_{n=1}^\infty \left(x_n,y\right)V^{-1}y_{n}. \eqno (11)$$
If we put $x = {\left(V^{-1}\right)}^{\ast}x$ in inequality (1), we obtain
$$A\left\| {\left(V^{-1}\right)}^{\ast}x \right\|_{X} \leq \left\| \Bigl(\left({\left(V^{-1}\right)}^{\ast}x, y_{n} \right)\Bigr) \right\|_{X_{d}} \leq B\left\| {\left(V^{-1}\right)}^{\ast}x \right\|_{X}. \eqno (12)$$
Since the operator $V^{-1}$ is invertible and bounded, (12) yields
$$\widetilde{A}\left\| x \right\|_{X} \leq A\left\| {\left(V^{-1}\right)}^{\ast}x \right\|_{X} \leq
\left\| \Bigl(\left({\left(V^{-1}\right)}^{\ast}x, y_{n} \right)\Bigr) \right\|_{X_{d}} \leq B\left\| {\left(V^{-1}\right)}^{\ast}x \right\|_{X}
\leq \widetilde{B}\left\| x \right\|_{X}.$$
And with the formula $\left(x, V^{-1}y_n\right) = \left({\left(V^{-1}\right)}^{\ast}x, y_n\right)$ we obtain inequality (1) for the sequence $\left(V^{-1}y_n\right)$.
To show that sequences $\left(x_n\right)$ and $\left(V^{-1}y_n\right)$ form a cross-frame, 
it remains to prove that from the formula (11) we can obtain the reconstruction formula
$$x = \sum_{n=1}^\infty \left(x,V^{-1}y_{n}\right)x_{n}. \eqno (13)$$
We define synthesis and analysis operators for the sequences $\left\{x_n\right\}$ and $\left\{ V^{-1}y_{n} \right\}$.
As in the Proposition 4, consider 
$$\widetilde{S}R = I_{X^\ast}. \eqno (14)$$
Applying the conjugation in (14), we obtain the corresponding equality for adjoint operators:
$$R^{\ast}{\widetilde{S}}^{\ast} = I_{{X^\ast}^\ast}.$$
From the above formula we obtain, that for all $x \in X$:
$$R^{\ast}\widetilde{R}x = x,$$
and with $\widetilde{R}x=\Bigl\{\left(x,V^{-1}y_{n}\right)\Bigr\}$, finally we get the desired 
reconstruction formula (13).\\
The theorem is proved.
}

\smallskip

Next, we consider yet another property of cross-frames in Banach space, and prove that these frames are not 
satisfying Kuznetsova theorem (see below).

\smallskip

{\bf Definition 6.}
{\hskip 0.2 cm
Family of operators $\left\{T_t\right\}_{t \in R}$ on Banach space X is called one-parameter group of operators 
if it satisfies the conditions: 
$T_{t+s} = T_t T_s$, $T_0 = I_X$.
$\left\{T_t\right\}_{t}$ is called strongly continuous if for any $x$ from $X$ the condition $t \to t_0$ yields ${T_t}x \to {T_{t_0}}x$. 
$\left\{T_t\right\}_{t}$ is called uniformely bounded if there exists positive constant $C$ such that for all $t \in R$ and $x \in X$ the next inequality holds: $\left\| {T_t}x \right\| \leq C\left\| x \right\|$.
}

\smallskip

The following theorem was proved by Kuznetsova in [5] for the case of unconditional basis in Banach space.

{\bf Theorem 8.}
{\hskip 0.2 cm \it 
Let $\left\{ e_n \right\}$ be an unconditional basis in Banach space $X$. Then there exists such one-parameter group of operators $\left\{T_t\right\}_{t \in R}$ that $T_t e_n = e^{i {\lambda}_n t} e_n$.
}

\smallskip

We are going to show that the cross-frame has not this property.\\
First let us prove two auxiliary statements.

\smallskip

{\bf Definition 7.}
{\hskip 0.2 cm  
Given operators $A:X_d \mapsto X_d$ and $B:X \mapsto X$.
Operator $S:X_d \mapsto X$ is called {\it intertwining operator} for the pair $\left[A, B\right]$, if the following diagram is  commutative:
$$
\begin{CD}
X_d @>A>> X_d \\
@VSVV @VVSV \\
X @>B>> X
\end{CD}
$$
i.e. if we have $SA=BS$.
}

\vskip 0.5 cm

Consider associated with the cross-frame $\left\{x_n\right\}$, $\left\{y_n\right\}$ synthesis operator $S:X_d \mapsto X$, defined by $$Sa=\sum a_{n} x_{n}.$$
Let $N$ be the space of coefficients of null-series in the sequence $\left\{x_n\right\}$, which coincides with the kernel of synthesis operator $N = Ker\left(S\right)$.
Along with the synthesis operator we consider the analysis operator $\widetilde{R}:X \mapsto X_d$ acting by the formula
$$\widetilde{R}x = \left\{\left(x, y_{n}\right)\right\}$$

\vskip 0.4 cm
{\bf Therem 9.}
{\hskip 0.2 cm \it 
Given operator $A:X_d \mapsto X_d$, and $S:X_d \mapsto X$ --- synthesis operator.
Then the existence of such an operator $B:X \mapsto X$, that $S$ is intertwining for the pair $\left[A, B\right]$, is equivalent to inclusion $A\left(N\right) \subset N$. And the operator $B$ is unique.
}
{
\\ Proof. 
Necessity. Let such an operator $B$ exists. Consider an arbitrary null-series $\sum a_n x_n$. 
For a sequence $a = \left\{a_n\right\} \in N$ we define the sequence $\left\{b_n\right\}$ by the formula $b = Aa$.
Then
$$\sum b_n x_n = S\left\{b_n\right\} = SA\left\{a_n\right\}= BS\left\{a_n\right\} = B0_X = 0,$$
i.e. $b \in N$. Inclusion $A\left(N\right) \subset N$ is verified.\\
Sufficiency. Now let $A\left(N\right) \subset N$. For a vector $x \in X$ we choose $a \in X_d$ so that
$x = Sa$, and let $Bx = SAa$. We are to verify the correctness of definition of the operator $B$, i.e. that it does not matter which $a \in X_d$ was chosen.
Indeed, if $x = Sa_1 = Sa_2$, then $a_1 - a_2 \in N$, and hence $A\left(a_1 - a_2\right) \in N$, from where $SA\left(a_1 - a_2\right) = 0$.
Thus $SAa_1 = SAa_2$. The correctness of definition of the operator $B$ is verified. Equality $SA = BS$ is valid by the construction. 
We verify uniqueness of the operator $B$. Suppose there is another operator, $\widetilde{B}$ having the same property. Then $BS = SA = \widetilde{B}S$.
By surjectivity of synthesis operator, we have that the operators $B$ and $\widetilde{B}$ coincide on the whole space $X$.
}
\smallskip

\vskip 0.4 cm
{\bf Theorem 10.}
{\hskip 0.2 cm \it 
Given operator $B:X \mapsto X$, $S$ --- synthesis operator.
Then there exists an operator $A:X_d \mapsto X_d$, that $S$ is an intertwining for the pair $\left[A, B\right]$, and the general form of such an operator is given by
$$A = \widetilde{R}BS + A_0, \eqno (15)$$
where $A_0$ --- an arbitrary operator satisfying the condition $A_0\left(X_d\right) \subset N$.
}
{
\\ Proof. 
Let for operator $A:X_d \mapsto X_d$ we have $$BS = SA.$$ We show that then it has the form (15).
Act with te operator $\widetilde{R}$ on the right:
$$BS\widetilde{R}=SA\widetilde{R}. \eqno (16)$$
It is known (see [1]) that for the operators $S$ and $\widetilde{R}$ associated with cross-frame we have $S\widetilde{R} = I_{X}$ and 
$\widetilde{R}S = I_{X_{d}} - P$, where $P$ is a projector onto the space $N$.
Rewrite (16) with view of the above:
$$B=SA\widetilde{R}. \eqno (17)$$
Now consider the construction of $\widetilde{R}BS$. In view of (17)
$$\widetilde{R}BS = \widetilde{R}SA\widetilde{R}S.$$
This implies:
$$\widetilde{R}BS = \left(I - P\right)A\left(I - P\right) = A - AP - PA + PAP.$$
By Theorem 9, $A\left(N\right) \subset N$. Then we set $A_0 = AP + PA - PAP$. It is obvious that 
$A_0\left(X_d\right) \subset N$. Finally, we find $A = RBS + A_0$.\\
Now we show that for any choice of the operator $A$ of the form (15) the operator $S$ is an intertwining for pair $\left[A, B\right]$.
To do this, just compute the $SA$: 
$$SA = S\left(\widetilde{R}BS + A_0\right) = S\widetilde{R}BS + SA_0 = BS + \theta = BS.$$
The theorem is proved.
}
\smallskip

\smallskip
{\bf Theorem 11.}
{\hskip 0.2 cm \it 
Given cross-frame $\left\{ x_n \right\}$, $\left\{ y_n \right\}$ in the Banach space $X$, which is not a basis for $X$.
Let $\lambda = \left\{\lambda_n\right\}$ be a separable sequence, i.e. there exists such positive real constant $\delta$ that for all $n \neq m$ we have $|\lambda_n - \lambda_m| > \delta$. Then for the sequence $\left\{ x_n \right\}$ there are no such one-parameter group of operators $\left\{T_t\right\}_{t \in R}$ that $T_t x_n = e^{i {\lambda}_n t} x_n$.
}
{
\\Proof.
Suppose that the sequence $\left\{ x_n \right\}$ satisfies the property from Theorem 8, i.e. there exists such one-parameter group of operators $\left\{T_t\right\}_{t \in R}$ that $T_t x_n = e^{i {\lambda}_n t} x_n$. 
Then by Theorem 10 for any $t \in R$ we can find such an operator $V_t : X_d \mapsto X_d$ that the synthesis operator
is intertwining for a pair $[V_t, T_t]$, and we know the view of this operator: $V_t = \widetilde{R} T_t S + A_{0t}$, where $A_{0t}$ is an arbitrary operator satisfying the condition $A_{0t}\left(X_d\right) \subset N$. We choose $A_{0t} = P e^{i \lambda_n t}$, where P is projector on $N$ by $Im \left( \widetilde{R} \right)$.\\
Let us compute the values of operators $V_t$ on basis vectors $\left\{ \varepsilon_n \right\}$: \\
$V_t \varepsilon_n = \widetilde{R} T_t S \varepsilon_n + P e^{i \lambda_n t} \varepsilon_n = 
 \widetilde{R} T_t x_n + P e^{i \lambda_n t} \varepsilon_n = 
 \widetilde{R} e^{i \lambda_n t} x_n + P e^{i \lambda_n t} \varepsilon_n = \\
 \widetilde{R} S e^{i \lambda_n t} \varepsilon_n + P e^{i \lambda_n t} \varepsilon_n = 
 \left(I_X - P\right) e^{i \lambda_n t} \varepsilon_n + P e^{i \lambda_n t} \varepsilon_n =
 e^{i \lambda_n t} \varepsilon_n$.\\
Consider an arbitrary null series $\sum a_n x_n$.
From Theorem 9 we know that $V_tN \subset N$, then $V_t {\sum a_n x_n} = 0$, i.e. 
$\sum e^{i \lambda_n t}a_n x_n = 0$. From the requirement of separability of $\left\{ \lambda_n \right\}$ 
with use of results regarding linear independence of exponential systems [6] we get that the sequence $\left\{ e^{i \lambda_n t} \right\}$ is w-linear independent.
Hence, all coefficients $a_n x_n = 0$, and since $x_n \neq 0$, we have $a_n = 0$. 
This means that the space of null series coefficients $N = \left\{ 0 \right\}$. But this happens only in the case 
when $\left\{ x_n \right\}$ forms a basis in $X$. The contradiction yields that the cross-frames in Banach space
do not satisfy the property from Theorem 8.
}
\smallskip

\vskip 1.2 cm
Bibliography:
\vskip 0.2 cm
\noindent 1. Casazza P., Han D., Larson D., Frames for Banach spaces, Contemp. Math., 247 (1999), 149-182.\\
\noindent 2. Grochenig K., Describing functions: atomic decompositions versus frames, Monat. Math., 112 (1991), 1-41.\\
\noindent 3. Casazza P., Christensen O., Stoeva D. T. Frame expansions in separable Banach spaces, J. Math. Anal. Appl., 307, 710-723 (2005).\\
\noindent 4. P. A. Terekhin, Frames in Banach Spaces, Funkts. Anal. Prilozh., 44:3 (2010), 50-62.\\
\noindent 5. Kuznetsova T.A., Bounded group of operators and approximation theory in the complex domains, Computational methods and programming: Interuniv. scient. collection of papers Saratov: Publishing of Saratov Univ (1981), 53-62.\\
\noindent 6. Leont'ev A.F., Entire functions: Series of exponentials, Moscow, (1983).
\end{document}